\documentstyle{amsppt}
\loadeusm
\topmatter
\title
A Probable Prime Test With High Confidence
\endtitle
\email grantham\@super.org
\endemail
\author Jon Grantham \endauthor
\endsubjclass
\thanks
This paper appeared in {\it Journal of Number Theory} {\bf 72}, 32--47
(1998).  Copyright \copyright 1998 by Academic Press.
\endthanks
\abstract
Monier and Rabin proved
that an odd composite can pass the Strong Probable Prime Test 
for at most $\frac 14$ of the possible bases.  In this paper,
a probable prime test is developed using quadratic polynomials and
the Frobenius automorphism.  The test, along with a fixed number of
trial divisions, ensures that a composite $n$ will pass for less than
$\frac 1{7710}$ of the polynomials $x^2-bx-c$ with
$\left(b^2+4c\over n\right)=-1$ and $\left(-c\over n\right)=1$.
The running time of the test is asymptotically
$3$ times that of the Strong Probable Prime Test.
\endabstract

\address Institute for Defense Analyses, Center for Computing Sciences,
17100 Science Drive, Bowie, MD 20715 
\endaddress
\NoBlackBoxes
\def\as{1}
\def\arnault{2}
\def\atkin{3}
\def\bbcgp{4}
\def\me{5}
\def\guy{6}
\def\knuth{7}
\def\mojones{8}
\def\mon{9}
\def\mont{10}
\def\dopo{11}
\def\psw{12}
\def\rabin{13}
\def\gcmd{\operatorname{gcmd}}

\endtopmatter
\document
\head{\S 1 Background}
\endhead
Perhaps the most common method for determining whether or not a number
is prime is the Strong Probable Prime Test.  Given an odd integer $n$, let
$n=2^rs+1$ with $s$ odd.  Choose a random integer $a$ with $1\le a\le n-1$.
If $a^s\equiv 1 \bmod n$ or $a^{2^js}\equiv -1 \bmod n$ for some $0\le
j\le r-1$, then $n$ passes the test.  An odd prime will pass
the test for all $a$.

The test is very fast; it requires no more than $(1+o(1))\log_2 n$
multiplications mod $n$, where $\log_2 n$ denotes the base $2$
logarithm.
 The catch is that a number which passes the
test is not necessarily prime.  Monier \cite{\mon} and Rabin
\cite{\rabin}, however, showed that a composite $n$ passes for at most
$\frac 14$ of the possible bases $a$.
Thus, if the bases $a$ are chosen
at random, composite $n$ will pass $k$ iterations of the Strong
Probable Prime Test with probability at most $\frac 1{4^k}$.

Recently, Arnault \cite{\arnault} has shown that any composite $n$
passes the Strong Lucas Probable Prime Test for at most $\frac 4{15}$ of the
bases $(b,c)$, unless $n$ is the product of twin primes having certain
properties (these composites are easy to detect).
Also, Jones and Mo \cite{\mojones} have introduced an Extra Strong
Lucas Probable Prime Test.  Composite $n$ will pass this test with
probability at most $\frac 18$ for a random choice of the parameter
$b$.  These authors did not concern themselves
with the issue of running
time. The methods of \cite{\mont}, however, show each test can be
performed in twice the time it takes to perform the Strong Probable
Prime Test.  By contrast, two iterations of the Strong Probable Prime Test
take the same amount of time and are passed by a composite
with probability at most $\frac 1{16}$.

Pomerance, Selfridge and Wagstaff have proposed a test, based on a
combination of the Strong Probable Prime Test and the Lucas Probable Prime
Test, that seems very powerful \cite{\psw}.
   Indeed, nobody has yet claimed the
\$620 that they offer for a composite that passes it \cite{\guy}, even
though they have relaxed the conditions of the offer so that the
Probable Prime Test no longer need be ``Strong''.
Adams and Shanks proposed a test
based on Perrin's sequence \cite{\as}.  
The Q and I cases of their test also have no known
pseudoprimes.   The success of both of these tests suggests that it is
possible to construct a test that is, in some sense, stronger than the
Strong Probable Prime Test.

My goal is to provide a test which is always passed 
by primes, but passed by
composites with probability less than $\frac 1{7710}$.  This test
has running time bounded by $(3+o(1))\log_2 n$ multiplications mod $n$.  By comparison, $3$
iterations of the Strong Probable Prime Test have a 
running time bounded by $(3+o(1))\log_2 n$
multiplications and have a probability of error at most $\frac
1{64}$.  

The number $\frac 1{7710}$ is slightly greater than $\frac{17}{2^{17}}+\epsilon$,
and $\frac{17}{2^{17}}$ is the limit of the analysis in this
paper.  The lack of counterexamples to the PSW and Perrin tests
indicates that the true probability of error may be much lower.
See Section 4 for possible improved versions of the test.

I would like to thank Carl Pomerance for his comments on this paper,
including the suggestion of the technique described in Proposition
3.3.

\head{\S 2 The Quadratic Frobenius Test}
\endhead
The following test is based on the concepts of Frobenius probable
primes
and strong Frobenius probable primes, which I introduced in another paper
\cite{\me}.
These tests involve computations in the ring $\Bbb{Z}[x]/(n,f(x))$,
where $f(x)\in\Bbb{Z}[x]$ and $n$ is an odd positive integer.  We will
be considering the special case $f(x)=x^2-bx-c$.
For convenience of notation, we introduce the integer $B$, which we later show 
can be taken to be $50000$.

\proclaim{Definition}
Suppose $n>1$ is odd, $\left(b^2+4c\over n\right)=-1$, and
$\left(-c\over n\right)=1$.
The {\bf Quadratic Frobenius Test (QFT) with parameters $(b,c)$} consists
of the following.

1) Test $n$ for divisibility by primes less than or equal to
$\min\{B,\sqrt{n}\}$.
If it is divisible by one of these primes, declare $n$ to be
composite and stop.
 
2) Test whether $\sqrt{n}\in\Bbb{Z}$.  If it is, declare $n$ to be composite
and stop.

3) Compute $x^{\frac{n+1}2} \bmod (n,x^2-bx-c)$.  If
$x^{\frac{n+1}2}\not\in\Bbb{Z}/n\Bbb{Z}$, declare $n$ to be composite and
stop.
 
4) Compute $\topsmash{x^{n+1}} \bmod (n,x^2-bx-c)$.  If 
$\topsmash{x^{n+1}}\not\equiv -c$,
declare $n$ to be composite and stop.

5) Let $\topsmash{n^2-1=2^rs}$, where $s$ is odd.
If $x^s\not\equiv 1 \bmod (n,x^2-bx-c)$,
and $\topsmash{x^{2^js}\not\equiv -1 \bmod (n,x^2-bx-c)}$ for all $0\le j\le
r-2$, declare $n$ to be composite and stop.

If $n$ is not declared composite in Steps 1--5, declare $n$ to be a
probable prime.
\endproclaim

\proclaim{Definition}
Suppose $n,b,c$ are as above.  We say $n$ passes the QFT with
parameters $(b,c)$ if the test declares $n$ to be a probable prime.
\endproclaim

The requirement that $\left(-c\over n\right)=1$ complicates the test
somewhat.  Without it and Step 3,
an accuracy of only $\frac 1{495}$ can be proven.

Steps 4 and 5 of the QFT are equivalent to the Strong Frobenius Probable
Prime Test described in \cite{\me} for the polynomial $x^2-bx-c$.
As is stated there, the test is not
entirely new, but rather a combination of ideas contained in earlier
probable prime tests.  The stricter condition $j\le r-2$ is possible because
the restriction $\left(-c\over n\right)=1$ forces $x$ to be a square in
the finite field $\Bbb{Z}[x]/(n,x^2-bx-c)$, if $n$ is a prime.

Traditionalists might prefer to rephrase the test in terms of Lucas
sequences.  For example, Step 4 is equivalent to $U_{n+1}(b,c)\equiv
0$ and $V_{n+1}(b,c)\equiv -2c \bmod n$, where $U_k$ and $V_k$ are the
standard Lucas sequences.
While this rephrasing might be useful to understand what is going on
in the test, the proof of accuracy relies on properties of finite
fields, so I do not feel it beneficial to use this notation.

\proclaim{Definition}
We define one iteration of the 
{\bf Random Quadratic Frobenius Test (RQFT)} for an odd integer $n>1$
to consist of the following:

1) Choose pairs $(b,c)$ at random with $1\le b,c < n$ until one is found with
$\left({b^2+4c}\over n\right)=-1$ and $\left(-c\over n\right)=1$, or with 
$\gcd(b^2+4c,n)$, $\gcd(b,n)$ or $\gcd(c,n)$ a nontrivial divisor of $n$.  
However, if 
the latter case occurs before the former, declare $n$ to be
composite and stop.
If after $B$ pairs are tested, none is found satisfying the
above conditions, declare $n$ to be a probable prime and stop.

2) Perform the QFT with parameters $(b,c)$.
\endproclaim

Of course, if more than one iteration of the RQFT is performed, Steps
1 and 2 of the QFT can be omitted in subsequent iterations.

Step 1 of the RQFT requires a declaration of probable primality if one is
extremely unlucky in choosing pairs $(b,c)$.  An objection could be
made that the declaration is not based on any actual
evidence that $n$ is prime.  This objection would be an accurate one,
but the declaration is only done in the extremely unlikely case that no
suitable pair is found.  Since we want to be certain to declare primes to be
probable primes, we must declare $n$ to be a probable prime.

Without this limit on the number of pairs, the running time of the test
is not deterministic.  Purists are welcome to delete this portion of
the test, but they will be left with a probabilistic running time.

\proclaim{Proposition 2.1}
If $p$ is an odd prime with $\topsmash{\left({b^2+4c}\over
p\right)}=-1$ and $\topsmash{\left(-c\over p\right)}=1$,
then $p$ passes the Quadratic Frobenius Test with
parameters $(b,c)$.
\endproclaim
\demo{Proof}
Proposition 2.1 follows from elementary properties of finite
fields.  Most of these can be found in \cite{\me}.  The only
differences here are the condition that $\left(-c\over p\right)=1$,
Step 3, and the elimination of the possibility that $x^{2^{r-1}s}\equiv -1$.

Since $x^{p+1}\equiv -c$ and $\left(-c \over p\right)=1$, we have
$x^{\frac{p+1}2} \in \Bbb{Z}/p\Bbb{Z}$, so $p$ passes Step 3.

It remains to show that $\left(-c\over p\right)=1$ implies that
$x^{2^{r-1}s}\equiv 1$.  Note that $x^{2^{r-1}s}=x^{\frac{p^2-1}2}
= (x^{(p+1)/2})^{p-1} \equiv 1$, since $x^{\frac{p+1}2} \in
\Bbb{Z}/p\Bbb{Z}$.
\enddemo

\proclaim{Proposition 2.2}
Let $p$ be an odd prime.
Let $\epsilon_1,\epsilon_2\in\{-1,1\}$.
If $\epsilon_1 \neq\epsilon_2$,
there are $\frac{(p-1)^2}4$ pairs $(b,c) \bmod p$ such that
$\botsmash{\left({b^2+4c}\over p\right)}=\epsilon_1$ 
and $\botsmash{\left(-c\over p\right)}=\epsilon_2$, otherwise there are
at most
$\frac{(p-1)^2}4-\epsilon_1\left(\frac{p-1}2\right)$ such pairs.
\endproclaim
\demo{Proof}
Pick any non-square $r \bmod p$.
Let $R=1$ if $\epsilon_2=1$, and $R=r$ if $\epsilon_2=-1$.  Let
$S=1$ if $\epsilon_1=1$ and $S=r$ if $\epsilon_1=-1$.
For each pair $(b,c)$ there exist four pairs $(x,y)$ 
with $-c\equiv Rx^2$ and $b^2+4c\equiv Sy^2$.  Substituting,
$b^2=R(2x)^2+Sy^2$.

We need to count
$$N=\frac 14\#\{(b,x,y)\bmod p : x,y\not\equiv 0\bmod p
\ \text{and}\ b^2 \equiv R(2x)^2+Sy^2\bmod p\}.$$

If $\epsilon_1=\epsilon_2=1$, then
$$\multline N=N_1=\frac 14\#\{(x,y,z)\ : x,y\not\equiv 0\bmod p \ \text{and}\ x^2+y^2 \equiv z^2\}\\
=\frac 14\#\{(x,y,z)\ : x,y\not\equiv 0\bmod p \ \text{and}\ (z-y)(z+y)\equiv x^2\}\\
=\frac 14\sum_{x=1}^{p-1}\#\{ab\equiv x^2 : a\not\equiv b\}
=(p-1)(p-1-2)/4=(p-1)(p-3)/4.  
\endmultline$$

If $\epsilon_1=\epsilon_2=-1$, then
$$\multline
N=N_2=\frac 14\#\{(x,y,z)\ : x,y\not\equiv 0\bmod p
\ \text{and}\ rx^2+ry^2 \equiv z^2\}\\
=\frac 14\#\{(x,y,z)\ : x,y\not\equiv 0\bmod p
\ \text{and}\ x^2+y^2 \equiv r(z/r)^2\}.
\endmultline$$
$$\multline
N_1+N_2
=\frac 14\#\{(x,y,z)\ : x,y\not\equiv 0\bmod p
\ \text{and}\ x^2+y^2\equiv z^2\text{ or }rz^2\}\\
\le\frac 12\#\{(x,y)\ : x,y\not\equiv 0\bmod p\}\le(p-1)^2/2.
\endmultline$$
So $N_2\le(p^2-1)/4+(p-1)/2$.

If $\epsilon_1\neq\epsilon_2$, then
$$\multline
N=N_3=\frac 14\#\{(x,y,z)\ : x,y\not\equiv 0\bmod p
\ \text{and}\ x^2+ry^2 \equiv z^2\}\\
=\frac 14\#\{(x,y,z)\ : x,y\not\equiv 0\
\bmod p \ \text{and}\ x^2+ry^2 \equiv rz^2\}\\
=\frac 14\#\{(x,y)\ :
x,y\not\equiv 0\bmod p\}=(p-1)^2/4.
\endmultline$$
\enddemo

\proclaim{Proposition 2.3}
Let $n$ be an odd squarefree number.
Let $\epsilon_1,\epsilon_2\in\{-1,1\}$.  
If $\epsilon_1\neq\epsilon_2$, there are at most
$\frac 14\phi(n)^2$ pairs $(b,c) \bmod n$ such that $\left(b^2+4c\over n\right)=\epsilon_1$ and 
$\left(-c\over n\right)=\epsilon_2$, otherwise there are at most $\frac
14\phi(n)^2+\frac{\epsilon_1} 2\mu(n)\phi(n)$ such pairs.
\endproclaim
\demo{Proof}
Let $n=p_1\dots p_k$, where the $p_i$ are primes.  The case $k=1$
follows from Proposition 2.2.

We proceed by induction on $k$.  Let $N_{\epsilon_1,\epsilon_2}(n)$ be
the number of pairs  $(b,c) \bmod n$ such
that $\left(b^2+4c\over n\right)= \epsilon_1$ and
$\left(-c\over n\right)=\epsilon_2$.

Let $n=mp$.  Then, by the Chinese Remainder Theorem,
$$\multline
N_{\epsilon_1,\epsilon_2}(n)=N_{\epsilon_1,\epsilon_2}(m)N_{1,1}(p)+
N_{-\epsilon_1,\epsilon_2}(m)N_{-1,1}(p)\\
+N_{\epsilon_1,-\epsilon_2}(m)N_{1,-1}(p)
+N_{-\epsilon_1,-\epsilon_2}(m)N_{-1,-1}(p).
\endmultline$$

By the inductive hypothesis, if $\epsilon_1=\epsilon_2$,
$$\multline
N_{\epsilon_1,\epsilon_2}(n)
\le\left[\frac 14\phi(m)^2+\frac{\epsilon_1} 2\mu(m)\phi(m)\right]
\left[\frac{(p-1)^2}4-\frac{p-1}2\right]
+\frac 1{16}\phi(m)^2(p-1)^2\\
+\frac 1{16}\phi(m)^2(p-1)^2
+\left[\frac 14\phi(m)^2-\frac{\epsilon_1} 2\mu(m)\phi(m)\right]
\left[\frac{(p-1)^2}4+\frac{p-1}2\right]\\
=\frac 14\phi(m)^2(p-1)^2
-\frac{\epsilon_1\mu(m)\phi(m)(p-1)}2
=\frac 14\phi(mp)^2+\frac{\epsilon_1}2\mu(mp)\phi(mp).
\endmultline$$

If $\epsilon_1\neq\epsilon_2$,
$$\multline
N_{\epsilon_1,\epsilon_2}(n)\le\frac 14\phi(m)^2
\topsmash{\left[\frac{(p-1)^2}4-\frac{p-1}2\right]
+\left[\frac 14\phi(m)^2-\frac{\epsilon_1} 2\mu(m)\phi(m)\right]\frac{(p-1)^2}4}\\
+\left[\frac 14\phi(m)^2+\frac{\epsilon_1} 2\mu(m)\phi(m)\right]\frac{(p-1)^2}4
+\frac 14\phi(m)^2\left[\frac{(p-1)^2}4+\frac{p-1}2\right]\\
=\frac 14\phi(m)^2(p-1)^2=\frac 14\phi(mp)^2.
\endmultline$$
\enddemo

\proclaim{Proposition 2.4}
Let $n$ be an odd composite, not a perfect square.  
Let $M(n)$ be the number of
pairs $(b,c) \bmod n$ such that $\left(b^2+4c\over n\right)=-1$ and
$\left(-c\over n\right)=1$, $n>(b^2+4c,n)>1$, or $n>(c,n)>1$.
Then $M(n)>\frac{n^2}4$.
\endproclaim
\demo{Proof}
If $n$ is squarefree, exactly 
$\frac 34\phi(n)^2$ pairs have $\left(b^2+4c\over n\right)$ and
$\left(-c\over n\right)$ non-zero, but not equal to the specified values.
Note that the number of such pairs mod $np^2$ is $\frac 34\phi(np^2)^2$,
since each such pair mod $n$ corresponds to $\frac{\phi(np^2)^2}{\phi(n)^2}$
such pairs mod $np^2$.  Thus we can remove the restriction that $n$ be
squarefree.

There are $n^2$ pairs $(b,c) \bmod n$.
Exactly $n$ pairs have $n|b^2+4c$ and exactly $n$ have $n|c$.
There is an overlap of one pair, $(0,0)$.

So $M(n)>n^2-\frac 34\phi(n)^2-2n>n^2-\frac 34(n-2)^2-2n
=\frac {n^2}4+n-3>\frac {n^2}4$.
\enddemo

\proclaim{Corollary 2.5}
The probability of failing to find a suitable pair $(b,c)$ in Step 1 of the
RQFT is less than $(3/4)^B$.
\endproclaim

\proclaim{Definition}
A positive integer $n$ is said to pass the Random Quadratic Frobenius Test
with probability $\alpha$ if the number of $(b,c)$ with $0\le b,c\le n$
such that $n$ passes the
Quadratic Frobenius Test with parameters $(b,c)$ is equal to
$\alpha M(n)$.
\endproclaim

\proclaim{Theorem 2.6}
An odd composite passes the RQFT with probability less
than $\frac 1{7710}$.
\endproclaim

The proof of Theorem 2.6 will be given in a sequence of lemmas.

\proclaim{Lemma 2.7}
Let $n$ be an odd integer.
If $p$ is a prime such that $p^2$ divides $n$, then $n$ passes the RQFT with
probability less than $\frac{4}p$.
\endproclaim
\demo{Proof}
Let $k$ be such that $\topsmash{p^k}|n$, but $\topsmash{p^{k+1}}\nmid n$.
If $n$ passes the QFT with parameters $(b,c)$, then
$x^{n+1}\equiv -c \bmod (n,x^2-bx-c)$.
By (5) in the definition of the QFT, $x^{n^2-1}\equiv 1 \bmod (n,x^2-bx-c)$.
So ${c}^{n-1}\equiv 
1 \bmod p^k$.  There can be at most $\gcd(n-1,p^k-p^{k-1})=\gcd(n-1,p-1)
\le{p-1}$ solutions to this
congruence mod $p^k$.  Hence there are at most $p-1$ choices for $c \bmod p^k$ for which it is possible that  $n$ passes the QFT with parameters $(b,c)$,
for some $b$.
Thus there are at most $p^{k+1}-p^k$ pairs $(b,c) \bmod p^k$ such that $n$
passes the QFT with parameters $(b,c)$.

By the Chinese Remainder Theorem, each pair mod $p^k$
corresponds to $(n/p^k)^2$ pairs
mod $n$.  Thus $n$ passes for
at most $\botsmash{(p^{k+1}-p^k)\frac{n^2}{p^{2k}}=(1-\frac 1p)\frac
{n^2}{p^{k-1}}<\frac 4{p^{k-1}}M(n)}$ pairs by
Proposition 2.4.

The lemma follows since $k\ge 2$.
\enddemo

\proclaim{Lemma 2.8}
Let $n$ be an odd composite with $p|n$.
We write that ``$n$ passes the QFT with parameters
$(b,c) \bmod p$'' if $n$ passes the QFT for some parameters $(b',c')$,
with $(b',c')\equiv (b,c) \bmod p$.
There are at most $\frac{p-1}2$ distinct pairs $(b,c) \bmod p$ with
$\left({b^2+4c}\over p\right)=1$ such that $n$ passes the
Quadratic Frobenius Test with
parameters $(b,c) \bmod p$.
\endproclaim
\demo{Proof}
First, we count the number of pairs $(b,c) \bmod p$ with
$\left({b^2+4c}\over p\right)=1$ such that $x^{n+1} \equiv -c \bmod
(p,x^2-bx-c)$.  If $a_1$ and $a_2$ are the two roots of $x^2-bx-c
\bmod p$, then the above congruence gives
$x^{n+1} \equiv -c \bmod (p,x-a_1)$.  This congruence
is equivalent to $a_1^{n+1}\equiv -c
\equiv a_1a_2 \bmod p$, or $a_1^n\equiv a_2$, since $p\nmid a_1$.
Similarly, $a_2^n\equiv a_1 \bmod p$.

Thus the number of polynomials $x^2-bx-c$ with $x^n\equiv b-x$ is
equal to the number of sets of integers $\{a_1,a_2\} \bmod p$ with
$a_1^n\equiv a_2 \bmod p$ and $a_2^n\equiv a_1 \bmod p$.  This is no
more than $(p-1)/2$.
\enddemo

\proclaim{Lemma 2.9}
The number of pairs $(b,c) \bmod n$ for which a squarefree integer
$n$ with $k$ prime factors passes the QFT with parameters $(b,c)$, such that
$\left({b^2+4c}\over p\right)=1$ for some $p|n$, is less than
$\frac{n\phi(n)}{2B}$ if $k$ is even and less than
$\frac{n\phi(n)}{B^2}$ if $k$ is odd.
\endproclaim
\demo{Proof}
Write $n=p_1p_2\dots p_k$.  Note $k>1$, since $k=1$ would imply that
$\left({b^2+4c}\over n\right)=-1$ and $\left({b^2+4c}\over
p_1\right)=1$.  

We consider separately the cases where $\left({b^2+4c}\over p_i\right)=1$ for
exactly one $i$, and for more than one $i$.

{\bf Case 1. }
Consider those pairs $(b,c)$ for which $n$ passes
the QFT with parameters $(b,c)$ such that $\left({b^2+4c}\over p_i\right)=1$
for exactly one $p_i|n$.    Note that such an $n$ can only
have $\left({b^2+4c}\over n\right)=-1$ when $k$ is even.

There are at most $\frac{p_i-1}2$ pairs $(b,c) \bmod p_i$ for which $n$
can pass the QFT,
 and for $j\neq i$, at most $\frac{p_j^2-p_j}2$ pairs mod $p_j$ with 
$\left({b^2+4c}\over p_j\right)=-1$.

So by the Chinese Remainder Theorem, there are at most
$\sum_{1\le i\le k} \frac {n\phi(n)}{2^kp_i}$ pairs $(b,c)$ for
which $n$ passes the QFT and
$\left({b^2+4c}\over p_i\right)=1$
for exactly one $p_i|n$.  Since any prime $p_i|n$ must be at least
$B$, the total number of pairs is at most $\frac{kn\phi(n)}{2^kB}$.

{\bf Case 2. }
For each pair $(b,c)$, let $L(b,c)$ be the $k$-tuple whose $i$th coordinate is
$\left({b^2+4c}\over p_i\right)$.

The $k$-tuples in $\{1,-1\}^k$ with at least two $1$s correspond
to pairs $(b,c)$ with
$\left({b^2+4c}\over p\right)=1$ for more than one prime $p|n$.
Let $V$ denote this set of $k$-tuples.

Given an element of $L\in V$,
let $p_i<p_j$ be the largest two primes with $\left({b^2+4c}\over
p\right)=1$.
By Lemma 2.8, the number of pairs $(b,c)$ for which $n$ passes the
QFT mod $p_i$ is at most $\frac{p_i-1}2$, and similarly mod $p_j$.
The number of pairs mod $p_k$ for $k\neq i,j$ is at most
$p_k(p_k-1)/2$.
By the Chinese Remainder Theorem, the
total number of pairs $(b,c)$ with $L(b,c)=L$
such that $n$ passes the QFT with
parameters $(b,c)$ is  at most
$[(p_i-1)/2][(p_j-1)/2]\prod_{k\neq i,j}  p_k(p_k-1)/2
 =\frac{n\phi(n)}{2^kp_ip_j}$.  The total number of $k$-tuples
$L\in V$ is less than $2^k$.

So the number of pairs $(b,c)$ such that $\smash{\left({b^2+4c}\over p\right)=1}$
for more than one prime dividing
$n$ is less than or equal to 
$$\sum_{L\in V} \frac{n\phi(n)}{2^kB^2}
<\frac{n\phi(n)}{B^2}.$$

We now combine the analysis of the two cases to complete the proof of
the Lemma.

For $k=2$, we have $\left({b^2+4c}\over p_1\right)=1$ and
$\left({b^2+4c}\over p_2\right)=-1$, or vice versa, so Case 2 is
impossible.  The proof easily follows by taking $k=2$ in Case 1.

For $k>2$ even, we have shown that the total number of pairs $(b,c)$
mod $n$ for which $n$ can pass the QFT is less than 
$n\phi(n)\left(\frac {k}{2^kB}+\frac 1{B^2}\right)$.  This is less
than $\frac {n\phi(n)}{2B}$.

For $k$ odd, Case 1 cannot occur, so the Lemma is proven.
\enddemo

\proclaim{Corollary 2.10}
A squarefree integer $n$ with an even number of prime factors passes
the RQFT
with probability less than $\frac 2B$.
\endproclaim
\demo{Proof}
Combine Proposition 2.4 with Lemma 2.9.
\enddemo

Lemma 2.11 is a more exact form of Proposition 6.2 of \cite{\me} for quadratic
polynomials.

\proclaim{Lemma 2.11}
If an odd squarefree number $n$ has $3$ prime factors, the probability that it passes the
RQFT is less than $\frac 4{B^2}+\frac{3(B^2+1)}{2(B^4-3B^2)}$.
\endproclaim
\demo{Proof}
Write $n=p_1p_2p_3$.  
By Lemma 2.9,
we know that $n$ passes for at most $\frac{4M(n)}{B^2}$ pairs 
with $\left({b^2+4c}\over n\right)=-1$ and $\left({b^2+4c}\over
p_i\right)=1$ for some $i$.
So it suffices to show that $n$ passes the QFT for at most $\frac{3(B^2+1)}
{2(B^4-3B^2)}\frac{n^2}4$ pairs with
$\left({b^2+4c}\over p_i\right)=-1$ for all $i$.

We know that $x^{n+1}\equiv -c \bmod
(p_i,x^2-bx-c)$, and $x^{p_i+1}\equiv -c$.  Since $c$ is invertible
mod $n$, $x$ is invertible mod $(n,x^2-bx-c)$, and $x^{n-p_i}\equiv 1$.
For each $p_i$, we need to know how many solutions
there are in $\Bbb{F}_{p_i^2}$ to $y^{n-p_i}\equiv 1$.
This congruence holds for exactly
$\gcd(n-p_i,p_i^2-1)$ elements.  Each pair $(b,c)$ corresponds to two
elements with minimal polynomial $x^2-bx-c$,
which either both do or both do not satisfy $y^{n-p_i}\equiv 1$.

Let $\gcd(n-p_i,p_i^2-1)=k_i$.  
By the preceding paragraph, there are at most $k_i/2$ pairs $(b,c)
\bmod p_i$ with $x^{n-p_i}\equiv 1 \bmod x^2-bx-c \bmod p_i$.
Then by the Chinese Remainder Theorem, $n$ passes the QFT for at most
$\botsmash{k_1k_2k_3/8}$ pairs with all Jacobi symbols equal to  $-1$.

Since $(p_1,p_1^2-1)=1$,
$k_1=\gcd(p_2p_3-1,p_1^2-1)$, and similarly for $k_2$ and
$k_3$. Let $j_i=(p_i^2-1)/k_i$.
Let $r_1=(p_2p_3-1)/k_1$, and define $r_2$ and $r_3$ analogously.

Then $r_1(p_1^2-1)=j_1(p_2p_3-1)$,
$r_2(p_2^2-1)=j_2(p_1p_3-1)$, and
$r_3(p_3^2-1)=j_3(p_1p_2-1)$.

Let $C=j_1j_2j_3$.  Assume that $r_1r_2r_3>C$.
Then $\frac{r_1r_2r_3}{j_1j_2j_3}\ge 1+\frac 1C$.

Multiplying the three equalities above, we have that
$$\frac{r_1r_2r_3}{j_1j_2j_3}=\frac{(p_2p_3-1)(p_1p_3-1)(p_1p_2-1)}
{(p_1^2-1)(p_2^2-1)(p_3^2-1)}.$$

Multiplying the right hand side out,
$$\frac{r_1r_2r_3}{j_1j_2j_3}=\frac{n^2-p_1^2p_2p_3-p_1p_2^2p_3-p_1p_2p_3^2
+p_1p_2+p_2p_3+p_1p_3-1}{n^2-p_1^2p_2^2-p_1^2p_3^2-p_2^2p_3^2
+p_1^2+p_2^2+p_3^2-1}.$$

Thus
$$\frac{r_1r_2r_3}{j_1j_2j_3}<\frac{n^2+p_1p_2+p_2p_3+p_1p_3}
{n^2-p_1^2p_2^2-p_1^2p_3^2-p_2^2p_3^2}.$$

Therefore
$$\frac{1+(p_1p_2+p_2p_3+p_1p_2)/n^2}
{1-(p_1^2p_2^2+p_1^2p_3^2+p_2^2p_3^2)/n^2}>1+\frac 1C.$$

We now use the facts that $n=p_1p_2p_3$ and $p_i>B$ to obtain
$\frac{1+3/B^4}{1-3/B^2}> 1+\frac 1C$.

Thus $C>\frac{B^4-3B^2}{3B^2+3}$.

If $r_1r_2r_3<C$, then $\frac{r_1r_2r_3}{j_1j_2j_3}\le 1-\frac 1C$.
From above,
$$\frac{r_1r_2r_3}{j_1j_2j_3}>\frac{n^2-p_1^2p_2p_3-p_1p_2^2p_3-p_1p_2p_3^2}
{n^2+p_1^2+p_2^2+p_3^2}.$$
Therefore
$$\frac{1-(p_1^2p_2p_3+p_1p_2^2p_3+p_1p_2p_3^2)/n^2}
{1+(p_1^2+p_2^2+p_3^2)/n^2}<1-\frac 1C.$$
We now use the facts that $n=p_1p_2p_3$ and $p_i>B$ to obtain
$\frac{1-3/B^2}{1+3/B^4}<1-\frac 1C$.  Thus 
$C>\frac{B^4+3}{3B^2+3}>\frac{B^4-3B^2}{3B^2+3}$. 

We have that $k_1k_2k_3/8<\frac{n^2}{8C}$.  The lemma is proven unless
$j_1j_2j_3=r_1r_2r_3$.
We will now show that this condition is impossible.

Once again, we multiply out the three equalities involving the $r_i$
and $j_i$ terms, but this time we can cancel $r_1r_2r_3$ with
$j_1j_2j_3$.  We get
$$(p_1^2-1)(p_2^2-1)(p_3^2-1)=(p_2p_3-1)(p_1p_3-1)(p_1p_2-1).$$

But $(p_1^2 -1)(p_2^2 - 1) < (p_1 p_2 - 1)^2$.  Multiplying by
similar inequalities for $p_1, p_3$ and $p_2, p_3$, and taking the square root,
we have a contradiction.

Therefore, the only possibility is $j_1j_2j_3\neq r_1r_2r_3$, and we
have shown that this assumption gives the probability stated in the theorem.
\enddemo

\proclaim{Lemma 2.12}  If $n$ is squarefree and has $k$ prime factors, 
where $k$ is odd, $n$ passes  the RQFT
with probability less than $\frac 1{2^{3k-2}}+\frac 1{2^{4k-3}}+\frac 4{B^2}$.
\endproclaim
\demo{Proof}
By Lemma 2.9, the 
number of pairs $(b,c)$ such that $\left({b^2+4c}\over p\right)=1$ 
for some prime $p|n$, and $n$ passes the QFT with parameters $(b,c)$,
is less than $\frac{n^2}{B^2}$.  

Now assume that $\left({b^2+4c}\over p\right)=-1$ for all $p|n$.

Write $n=p_1p_2\dots p_k$.  Let $J$ be the largest integer 
such that $2^{J+1}|\gcd(p_1^2-1,\dots ,p_k^2-1)$.  Then for each $i$,
$p_i^2\equiv 1 \bmod 2^{J+1}$, so $n^2\equiv 1 \bmod 2^{J+1}$.  Thus
$J<r$.

The number of solutions that $y^{2^js}\equiv -1$ can have in
$\Bbb{F}_{p^2}$
is either $0$ or $\gcd(p^2-1,2^js)$.
Furthermore, it only has solutions if
$-1$ is a perfect $2^j$th power.   This will only happen if
$2^{j+1}|{p^2-1}$.

So there are no pairs $(b,c)$ for which $\topsmash{x^{2^js}\equiv -1 \bmod (n,x^2-bx-c)}$
if $j>J$.

For at least one prime $p|n$, $2^{J+1}$ is the highest power of $2$
dividing $p^2-1$.  
In Step 4 of the QFT, we show that $-c$ is a square
mod $n$, and therefore mod $p$.  By Proposition 2.1, it
follows that $x^{\frac{p^2-1}2}\equiv 1
\bmod (p,x^2-bx-c)$. So $x$ has order dividing $2^J\ell$, for some odd
number $\ell$.  If $x^{2^Js}\equiv -1 \bmod (p,x^2-bx-c)$, then $2^{J+1}$
divides the order of $x$, and we have a contradiction.  Thus there are
no pairs for which $x^{2^Js}\equiv -1 \bmod (n,x^2-bx-c)$.

We consider $0\le j<J$.
For each $p_i$, there are $\gcd(p_i^2-1,2^js)$ solutions to 
$y^{2^js}\equiv -1$ in $\Bbb{F}_{p_i^2}$.  Each solution 
$y\not\in\Bbb{F}_{p_i}$ produces a solution to $x^{2^js}\equiv -1 \bmod 
(p_i,x^2-bx-c)$, where $x^2-bx-c$ is the minimum polynomial of $y$.
But for each minimum polynomial, both of its roots in
$\Bbb{F}_{p_i^2}$ will be solutions.  (The roots are distinct by the
restriction on $c$.)

Therefore, there are at most $\gcd(p_i^2-1,2^js)/2$ pairs $(b,c)$ mod
$p_i$ for which
$x^{2^js}\equiv -1 \bmod (p_i,x^2-bx-c)$.  Let $G=\prod_{p_i}
\gcd(p_i^2-1,s)$.  By the Chinese Remainder
Theorem, there are at most $\prod_{1\le i\le k} \gcd(p_i^2-1,2^js)/2$
pairs mod $n$ for which $x^{2^js}\equiv -1$  mod 
$(n,x^2-bx-c)$. But $\prod_{1\le
i\le k} \gcd(p_i^2-1,2^js)/2=\prod_{1\le i\le k}2^{j-1}\gcd(p_i^2-1,s)=
2^{(j-1)k}G$.
Similarly,
the number of pairs $(b,c)$ for which $x^s\equiv 1$ mod $(n,x^2-bx-c)$
is at most $G/2^k$.
Thus the total
number of pairs for which $n$ passes the QFT is bounded by $G/2^k+
\sum_{j=0}^{J-1}
2^{(j-1)k}G=\left[1+\frac{2^{Jk}-1}{2^k-1}\right]\frac{G}{2^k}$.

Since $\gcd(p_i^2-1,s)<p_i^2/2^{J+1}$, the number of pairs for which
$n$ passes the QFT is less
than
$\left[1+\frac{2^{Jk}-1}{2^k-1}\right]\frac{n^2}{2^{(J+1)k}2^k}$.
We have $J\ge 2$, and this expression
 is maximized at $J=2$, where it is equal to
$(2^k+2)\frac{n^2}{2^{4k}}=n^2(\frac 1{2^{3k}}+\frac 1{2^{4k-1}})$.

The lemma now follows from Proposition 2.4.
\enddemo

\demo{Proof of Theorem 2.6}
Choose $B=50000$.  If $n$ is not squarefree, or has an even number of prime
factors, Lemma 2.7 and Corollary 2.10 prove Theorem 2.6.
If $n$ is squarefree and has
$3$ prime factors, we apply Lemma 2.11.  If $n$ is squarefree and
 has $k$ prime factors, $k>3$, we apply Lemma 2.12.  The bound for
the probability given by that lemma is largest when $k=5$ and is
$1/2^{13}+1/2^{17}+1/25000^2$.

Note that adding on the number $\topsmash{\left(\frac 34\right)^B}$ from 
Corollary 2.5 does not increase the probability above $\frac
1{7710}$.
\enddemo

\head{\S 3 Running Time}
\endhead
Atkin \cite{\atkin} has
suggested a unit of running time of probable prime tests
based on the running time of the Strong Probable Prime Test.  The unit
is named the ``selfridge'' to honor John Selfridge for his discovery
of the Strong Probable Prime Test.  We give a slightly different
formalization below.

\proclaim{Definition}
An algorithm with input $n$
 is said to have running time of $k$ {\bf selfridges} if it can be completed
in the time it takes to perform $(k+o(1))\log_2 n$ multiplications mod
$n$.  Calculation
of an inverse will take time equal to $O(1)$ multiplications and
computation of $\left({\ }\over n\right)$ will take $O(1)$ multiplications.
Here, $O(1)$ is bounded as $n\rightarrow\infty$.
We will assume that addition takes $o(1)$ multiplications, where
$o(1)\rightarrow 0$ as $n\rightarrow\infty$.
\endproclaim

If a multiplication algorithm that took $O(\log n)$ bit operations
were discovered, the last assumption would be invalidated.
The definition of ``selfridge,'' however, is meant to be
a practical, if not completely precise, method of comparing different probable
prime tests.  Since the development of a multiplication algorithm that
takes $O(\log n)$ bit operations appears to be highly
improbable, I make the definition with a clear conscience.

Note that we use the term ``multiplication'' to refer both to the
arithmetic operation and the time required to perform it.
The above
definition of the term ``selfridge'' is motivated by the following proposition.

\proclaim{Proposition 3.1}
The Strong Probable Prime Test has a running time of at most $1$ selfridge.
\endproclaim
\demo{Proof}
If we write $n=2^rs+1$, with $s$ odd, then the Strong Probable Prime
Test requires, at most, the raising of an integer to the $s$ power mod
$n$ and the completion of $r$ squarings (which are multiplications).
By \cite{\knuth}, exponentiation to
the $t$th power can be done in $(1+o(1))\log_2 t$ multiplications 
using easily constructed
addition chains.  Since $\log_2 n>\log_2 s+r$, the test can be performed
in $(1+o(1))\log_2 n$ multiplications mod $n$.  Note that for most odd $n$, the
running time of the test is equal to $1$ selfridge.
\enddemo

It is the goal of this section to show that the Quadratic Frobenius
Test has 
running time of at most $3$ selfridges.  One ``trick'' used to attain
this running time relies on the use of Lemma 4.8 of \cite{\me}.

\proclaim{Lemma 4.8 of \cite{\me}}
Let $m,n$ be positive integers, and let $f(x),g(x),r(x)\in\Bbb{Z}[x]$.
If $f(r(x))\equiv 0 \bmod (n,f(x))$ and
$x^m \equiv g(x) \bmod (n,f(x))$, then
$r(x)^m \equiv g(r(x)) \bmod (n,f(x))$.
\endproclaim

Before analyzing the running time, we need to know how many
multiplications mod $n$ it takes to perform a multiplication mod
$(n,x^2-bx-c)$.

\proclaim{Proposition 3.2}
Multiplication $\bmod (n,x^2-bx-c)$ can be done in at most $5+o(1)$
multiplications $\bmod\ n$.
\endproclaim
\demo{Proof}
To multiply $(dx+e)(fx+g)\equiv (dg+ef+bdf)x+eg+cdf \bmod (n,x^2-bx-c)$, we
compute the $5$ products $df$, $eg$, $b(df)$, $c(df)$, and $(d+e)(f+g)$.  Then
we can compute $dg+ef+bdf=(d+e)(f+g)+bdf-df-eg$ and $eg+cdf$ by addition and
subtraction.
\enddemo

\proclaim{Proposition 3.3}
Let $A_j=x^j+(b-x)^j \bmod (n,x^2-bx-c)$, and
let $B_j=\frac{x^j-(b-x)^j}{2x-b}$  $\bmod (n,x^2-bx-c)$.  Let $C_j=c^j
\bmod n$.
Given values of $2^{-1}$ and  $(b^2+4c)$, 
$\bmod (n,x^2-bx-c)$,
then $(A_{j+k},B_{j+k},C_{j+k})$ can be computed from $(A_j,B_j,C_j)$ 
and $(A_k,B_k,C_k)$
in $8+o(1)$ multiplications $\bmod\ n$. (We call this type of computation a
{\bf chain addition}.)  $(A_{2j},B_{2j},C_{2j})$ can be computed
from $(A_j,B_j,C_j)$ in $3+o(1)$ multiplications $\bmod\ n$.  (We call this
type of computation a {\bf doubling}.)
$(A_j,B_j,C_j)$ can be computed in $(3+o(1))\log_2 j$ multiplications
$\bmod n$.  $x^j$ can be
computed from $(A_j,B_j)$ in $2+o(1)$ multiplications mod $n$.
\endproclaim
\demo{Proof}
We have
$A_j,B_j\in\Bbb{Z}/n\Bbb{Z}$.  Note the identities
$A_{j+k}=2^{-1}(A_jA_k+(b^2+4c)B_jB_k)$, $B_{j+k}=2^{-1}(A_jB_k+A_kB_j)$,
and $C_{j+k}=C_jC_k$.

This shows a chain addition takes $8+o(1)$ multiplications mod $n$.

Also note that $A_{2j}=A_j^2-2(-1)^jC_j$, $B_{2j}=A_jB_j$, and $C_{2j}=C_j^2$.
$2(-1)^jC_j=C_j+C_j$ or $-C_j-C_j$ (depending on the parity of $j$), 
so this part of the computation is
$o(1)$ multiplications.  So a doubling can be achieved in $3+o(1)$
multiplications mod $n$.

Once we have $A_j$ and $B_j$, we can compute
$x^j=B_jx+2^{-1}(A_j-bB_j)$
with one multiplication of $b$ and $B_j$, one subtraction, and one
multiplication by $2^{-1}$.
Thus finding $x^j$ from the pair $(A_j,B_j)$ costs $2+o(1)$
multiplications mod $n$.

We can compute $(A_j,B_j,C_j)$ using $(1+o(1))\log_2 j$ steps (doublings
or chain additions) by the addition chain methods described in \cite{\knuth}.
Since $o(\log n)$ of these steps will not be doublings, we can compute
$(A_j,B_j,C_j)$ in $(3+o(1))\log_2 j$ multiplications mod $n$.
\enddemo

\proclaim{Theorem 3.4}
The Random Quadratic Frobenius Test has running time of at
most $3$
selfridges.
\endproclaim
\demo{Proof}
For a given $n$, Step 1 of the RQFT will take at most $B$
tries in searching for a suitable pair $(b,c)$.
Thus it takes $O(1)$
multiplications mod $n$.  It remains to show that the QFT has running
time of $3$ selfridges.

Steps 1 and 2 of the QFT have running time bounded by a fixed number of
multiplications mod $n$.  

Assume $n\equiv 1 \bmod 4$.
Write $n-1=2^{r'}s'$, with $s'$ odd.  Then $n^2-1=2^{r'+1}(2^{r'-1}{s'}^2+s')$.
So $r=r'+1$ and $s=2^{r'-1}{s'}^2+s'$.  Write $t=\frac{s'-1}2$.

Computing $(A_t,B_t,C_t)$ takes $(3+o(1))\log_2 t$ steps,
by Proposition 3.3.  
Computing $x^n$ from
$(A_t,B_t,C_t)$ can be accomplished by first computing $(A_{s'},B_{s'},C_{s'})$, which
requires $11$ multiplications mod $n$.
 We then 
perform $r'-1$ doublings to get 
$(A_{\frac{n-1}2},B_{\frac{n-1}2},C_{\frac{n-1}2})$.  We can compute
$x^{\frac{n-1}2}$ in $2+o(1)$ multiplications
and then $x^{\frac{n+1}2}$ in $5+o(1)$ more.  This calculation
completes Step 3.  Squaring this result, we complete Step 4 with
at most $5+o(1)$ additional multiplications.

So the total
number of multiplications mod $n$ is $(3+o(1))\log_2
t+3r'+15+o(1)=(3+o(1))\log_2 n$.  
So Steps 3 and 4 take $3$ selfridges when $n\equiv 1 \bmod 4$.

Note that Step 5 only needs to be performed if $x^{n+1}\equiv -c \bmod
(n,x^2-bx-c)$, which implies $x^n\equiv b-x$.

For Step 5, note that $s=nt+t+\frac{n-1}2+1$ and, for $0\le e<r$,
$2^{e+1}s=n2^es'+2^es'$.  

Let $\sigma$ be the map from $(\Bbb{Z}/n\Bbb{Z})[x]/(n,x^2-bx-c)$ to itself that
sends $x$ to $b-x$.
We have $x^{nt}\equiv (b-x)^t \equiv \sigma(x^t)$, by Lemma 4.8 of
\cite{\me}. 

So $x^s\equiv x^{nt}x^tx^{\frac{n-1}2}x
\equiv \sigma(x^t)x^tx^{\frac{n-1}2}x$.  We have already
computed $(A_t,B_t)$ and $x^{\frac{n-1}2}$
in the process of computing
$x^n$, so it takes $2+o(1)$ multiplications to find $x^t$ and
$x^{\frac{n-1}2}$.
 Given $x^t$, we require at most $1$ multiplication mod $n$ and
an addition to compute $\sigma(x^t)$.  It takes $3+o(1)$ multiplications mod
$(n,x^2-bx-c)$, or $15+o(1) \bmod n$ to multiply the results together to get
$x^s$.  So we can compute $x^s$ with at
most $18+o(1)$ additional multiplications mod $n$.

We have $x^{2^{e+1}s}\equiv x^{n2^es'}x^{2^es'}$.  By Lemma 4.8 of
\cite{\me}, 
$x^{n2^es'}\equiv \sigma(x^{2^es'})$.  We have already computed
$(A_{2^es'},B_{2^es'})$ 
in the process of computing $x^{n+1}$.  It takes $2+o(1)$ multiplications to
compute $x^{2^es'}$ from these.  We require $1$
multiplication mod $n$ and one addition
to apply $\sigma$, and $5+o(1)$ more to find the
product $x^{n2^es'}x^{2^es'}$.  So we need 
$8+o(1)$ multiplications mod $n$ to compute $x^{2^{e+1}s}$.

We want to know if $x^s\equiv\pm 1$ or if $x^{2^{e+1}s}\equiv -1$ for
some $e<r$.  This can be accomplished by a binary search.

We first check whether or not $\topsmash{x^s}\equiv\pm 1$.  This step takes $22+o(1)$ multiplications mod $n$.

We begin by computing $\topsmash{x^{2^js}} \bmod (n,x^2-bx-c)$ for $j=[\frac
r2]$.  If $x^{2^js}\equiv -1$, we can stop.  If $x^{2^js}\equiv 1$, we
know that $x^{2^{e+1}s}\equiv 1$ for $e\ge j$, so we only need to
check those $e$ less than $j$.  If $x^{2^js}\not\equiv 1$, then we
know that $x^{2^{e+1}s}\not\equiv -1$ for $e<j$, so we only need to
check those $e$ greater than or equal to  $j$.  We continue in a similar manner
until we either find $e$ with $x^{2^{e+1}s}\equiv -1$ or prove that
none exists.  This takes at most $\log_2 r+1$ steps, each of which
requires at most $8+o(1)$ multiplications.  $r<\log_2 n+1$, so the total
number of multiplications mod $n$ in Step 5 is $O(\log\log n)$ when
$n\equiv 1 \bmod 4$.

If $n\equiv -1 \bmod 4$, write $\smash{n+1=2^{r'}s'}$.  Then
$\smash{n^2-1=2^{r'+1}(2^{r'-1}{s'}^2-s')}$.  Write $t=\frac{s'-1}2$.

In Step 3, we compute $(A_t,B_t,C_t)$, and then 
$(A_{s'},B_{s'},C_{s'})$.  We double $r'-1$ times, and then compute
$x^{\frac{n+1}2}$.  We then square to get $x^{n+1}$.
By a similar
analysis to the case when $n\equiv 1 \bmod 4$, Steps 3 and 4
require $3$ selfridges when $n\equiv -1 \bmod 4$.

For Step 5, observe that $s=nt-t+\frac{n+1}2-1$ and
$2^es=n2^es'-{2^e}s'$.  We can calculate $x^t$ and
$x^\frac{n+1}2$, and we can use the fact that $x^{nt}\equiv \sigma(x^t)$.

We can then proceed with Step 5 via the binary search described in the
case where $n\equiv 1 \bmod 4$.  Again, Step 5 takes $O(\log\log n)$
multiplications mod $n$.
\enddemo

\head{\S 4 Unanswered Questions}
\endhead
The most obvious question that arises is, ``Can we do significantly
better than $7710$?''
Increasing $B$ will significantly improve each lemma except for Lemma
2.12.  In order to do significantly better than $7710$, then, it seems
that we need to develop a
better analysis of the case where $n$ has $5$ prime factors.  An
improved analysis of the case where $n$ has $k$ prime factors, for any
odd $k$, would be even better.

The answer to this question is probably ``yes''.  Part of the motivation for
this test was the combined Strong Probable Prime and Lucas Probable
Prime 
tests of Pomerance, Selfridge, and Wagstaff \cite{\psw}.  Integers that
have $\frac{p^2-1}k|n^2-1$, with $k$ small, for each $p|n$ seem to offer
the best chances of ``fooling'' either test, but such numbers have
proven difficult to construct in practice.  This statement does not
preclude the existence of many such numbers, but it is encouraging to
note that the best heuristics for constructing such numbers
\cite{\dopo} would produce numbers with many prime factors.

The following theorem shows that composites that pass the
test with high probability have special form.

\proclaim{Theorem 4.1}
Let $n$ be an odd composite with $k$ prime factors.  Let $r,s$ be such
that $n^2-1=2^rs$ with $s$ odd.  Let $G=\prod_{p|n} \gcd(p^2-1,s)$.
Let $J$ be the largest integer such that $2^{J+1}|p^2-1$ for every
$p|n$.
If $\frac{n^2}{2^{k(J+1)}G}>\frac B{2^{3k+1}}$, then $n$ passes the RQFT
with probability less than $\frac 4B$, for any $B\ge 50000$.
\endproclaim
\demo{Proof}
The theorem follows immediately unless $k$ is odd and greater than
$3$.  In the proof of Lemma 2.12 we used the bound
$G<\frac{n^2}{2^{k(J+1)}}$.  Here we have the bound
$G<\frac{2^{3k+1}}B\frac{n^2}{2^{k(J+1)}}$.  The multiple of 
$\frac{2^{3k+1}}B$ carries through to the final probability bound,
giving a probability of error less than $\frac{4}B$.
\enddemo

We showed that $\frac{n^2}{2^{k(J+1)}G}>1$.  It would be interesting
to see computational results for small $n$, and further theoretical results
that might improve Theorem 2.6.

Scott Contini asks if a result analogous to that in \cite{\bbcgp}  can be proven
about the reliability of the RQFT in generating primes.

Another question is whether it would be useful to construct cubic and
higher order versions of the Random Quadratic Frobenius Test.  A basis
for doing so can be found in \cite{\me}.
The fraction $\frac
1{7710}$ would have to be improved considerably to make the increased
running time in larger finite fields
worthwhile.  One technique for doing so would be to
exploit the fact that if $(n,m)=1$, then for $m$ with $\phi(m)|d$,
$m|n^d-1$.  We could actually improve the Quadratic Frobenius Test
by observing that $3|n^2-1$ for all $n$ coprime to $3$.  Checking
$\gcmd(x^{(n^2-1)/3}-1,x^2-bx-c)$ may expose $n$ as composite.  This
modification would,
however, further complicate the analysis of the running time, so it will be 
taken up in a future paper.

\enddocument